\documentclass[12pt]{amsart}
\usepackage{latexsym,amsmath,amssymb}

\title[Graphs in the Heisenberg group]{The Lusin theorem and horizontal graphs in the Heisenberg group}

\author{Piotr Haj\l{}asz, Jacob Mirra}

\address{Department of Mathematics, University of Pittsburgh, 301
  Thackeray Hall, Pittsburgh, PA 15260, USA, {\tt hajlasz@pitt.edu}}

\address{Department of Mathematics, University of Pittsburgh, 301
  Thackeray Hall, Pittsburgh, PA 15260, USA, {\tt jrm152@pitt.edu}}

\thanks{P.H. was supported by NSF grant DMS-1161425.}

plus 6pt minus 12pt
\def\eps{\varepsilon}

\newcommand{\boldg}{{\mathbf g}}

\def\H{{\mathcal H}}

\newtheorem{theorem}{Theorem}
\newtheorem{lemma}[theorem]{Lemma}

\newtheorem{proposition}[theorem]{Proposition}

\theoremstyle{definition}

\newcommand{\barint}{
\rule[.036in]{.12in}{.009in}\kern-.16in \displaystyle\int }

\newcommand{\barcal}{\mbox{$ \rule[.036in]{.11in}{.007in}\kern-.128in\int $}}

\newcommand{\bbbr}{\mathbb R}

\newcommand{\bbbh}{\mathbb H}

\def\bbbc{{\mathchoice {\setbox0=\hbox{$\displaystyle\rm C$}\hbox{\hbox
to0pt{\kern0.4\wd0\vrule height0.9\ht0\hss}\box0}}
{\setbox0=\hbox{$\textstyle\rm C$}\hbox{\hbox
to0pt{\kern0.4\wd0\vrule height0.9\ht0\hss}\box0}}
{\setbox0=\hbox{$\scriptstyle\rm C$}\hbox{\hbox
to0pt{\kern0.4\wd0\vrule height0.9\ht0\hss}\box0}}
{\setbox0=\hbox{$\scriptscriptstyle\rm C$}\hbox{\hbox
to0pt{\kern0.4\wd0\vrule height0.9\ht0\hss}\box0}}}}

\def\bbbq{{\mathchoice {\setbox0=\hbox{$\displaystyle\rm Q$}\hbox{\raise
0.15\ht0\hbox to0pt{\kern0.4\wd0\vrule height0.8\ht0\hss}\box0}}
{\setbox0=\hbox{$\textstyle\rm Q$}\hbox{\raise 0.15\ht0\hbox
to0pt{\kern0.4\wd0\vrule height0.8\ht0\hss}\box0}}
{\setbox0=\hbox{$\scriptstyle\rm Q$}\hbox{\raise 0.15\ht0\hbox
to0pt{\kern0.4\wd0\vrule height0.7\ht0\hss}\box0}}
{\setbox0=\hbox{$\scriptscriptstyle\rm Q$}\hbox{\raise
0.15\ht0\hbox to0pt{\kern0.4\wd0\vrule height0.7\ht0\hss}\box0}}}}

\def\bbbz{{\mathchoice {\hbox{$\sf\textstyle Z\kern-0.4em Z$}}
{\hbox{$\sf\textstyle Z\kern-0.4em Z$}} {\hbox{$\sf\scriptstyle
Z\kern-0.3em Z$}} {\hbox{$\sf\scriptscriptstyle Z\kern-0.2em
Z$}}}}

\def\mvint_#1{\mathchoice
          {\mathop{\vrule width 6pt height 3 pt depth -2.5pt
                  \kern -8pt \intop}\nolimits_{\kern -3pt #1}}%
          {\mathop{\vrule width 5pt height 3 pt depth -2.6pt
                  \kern -6pt \intop}\nolimits_{#1}}%
          {\mathop{\vrule width 5pt height 3 pt depth -2.6pt
                  \kern -6pt \intop}\nolimits_{#1}}%
          {\mathop{\vrule width 5pt height 3 pt depth -2.6pt
                  \kern -6pt \intop}\nolimits_{#1}}}

\numberwithin{theorem}{section} \numberwithin{equation}{section}

\begin{document}


\subjclass[2000]{Primary 46E35; Secondary 46E30}
\keywords{Lusin theorem, Heisenberg group, characteristic points}

\sloppy

\begin{abstract}
In this paper we prove that every collection of measurable
functions $f_\alpha$, $|\alpha|=m$ coincides a.e. with
$m$th order derivatives of a function $g\in C^{m-1}$ whose
derivatives of order $m-1$ may have
any modulus of continuity weaker than that of a Lipschitz function.
This is a stronger version of earlier results of Lusin, Moonens-Pfeffer and Francos.
As an application we construct surfaces in the Heisenberg group 
with tangent spaces being horizontal a.e.
\end{abstract}

\maketitle

\section{Introduction}

In 1917 Lusin \cite{lusin} proved that for every measurable function $f:\bbbr\to\bbbr$ there is 
a continuous function $g:\bbbr\to\bbbr$ that is differentiable a.e. and such that
$g'(x)=f(x)$ for almost all $x\in\bbbr$. This result was generalized by Moonens and Pfeffer
\cite{moonensp} to the case of functions defined in $\bbbr^n$ and then by Francos \cite{francos}
to the case of higher order derivatives in $\bbbr^n$. He proved that if $f_\alpha$, $|\alpha|=m$
are measurable functions in an open set $\Omega\subset\bbbr^n$, then there is a function 
$g\in C^{m-1}(\Omega)$ that is $m$ times differentiable a.e. and such that for all $|\alpha|=m$,
$D^\alpha g=f_\alpha$ a.e. It is easy to see that in general one cannot require that
$g\in C^{m-1,1}_{\rm loc}$, i.e. one cannot assume that the derivatives of order $m-1$ are
Lipschitz continuous. For example in the case $m=1$ one cannot find a 
locally Lipschitz continuous function $g$ on $\bbbr^2$ such that $\nabla g(x,y)=\langle 2y,-2x\rangle$. Indeed,
for such a function we would have
$$
g(1,1)=g(0,0)+\int_0^1\frac{\partial g}{\partial x}(t,0)\, dt + 
\int_0^1\frac{\partial g}{\partial y}(1,t)\, dt = g(0,0)-2,
$$
$$
g(1,1)=g(0,0)+\int_0^1\frac{\partial g}{\partial y}(0,t)\, dt + 
\int_0^1\frac{\partial g}{\partial x}(t,1)\, dt = g(0,0)+2
$$
which is impossible. 

Clearly, continuity of derivatives of order $m-1$ 
in Francos' theorem result from some
uniform convergence and one could expect that with keeping track of estimates it should be possible
to prove H\"older continuity of derivatives of order $m-1$. However, as we will see, a much stronger result
is true. Namely we shall prove that
it is possible to construct a function $g$
with {\em any} modulus of continuity 
of derivatives of order $m-1$
which is worse than that of a Lipschitz function. 
\begin{theorem}
\label{main}
Let $\Omega\subset\bbbr^n$ be open, $m\geq 1$ an integer,
and let $f=(f_{\alpha})_{|\alpha|=m}$, be a collection of measurable functions
$f_{\alpha}:\Omega\rightarrow\bbbr$, $|\alpha|=m$.
Let $\sigma>0$ and let
$\mu:[0,\infty)\rightarrow[0,\infty)$ be a continuous function
with $\mu(0)=0$ and $\mu(t)=O(t)$ as $t\rightarrow\infty$. Then
there is a function $g\in C^{m-1}(\bbbr^n)$ that is $m$-times differentiable
a.e., and such that
\begin{enumerate}
\item[(i)] $D^\alpha g=f_\alpha$ a.e. on $\Omega$ for all $|\alpha|=m$;
\item[(ii)] $\left\Vert D^\gamma g\right\Vert _{L^\infty(\bbbr^n)}<\sigma$ for all $0\leq |\gamma|\leq m-1$;
\item[(iii)] 
$$
|D^\gamma g(x)-D^\gamma g(y)|\leq \sigma |x-y|
$$
for all $x,y\in\bbbr^n$ and all $0\leq |\gamma|\leq m-2$;
\item[(iv)] 
$$
|D^\gamma g(x)-D^\gamma g(y)|\leq\frac{|x-y|}{\mu(|x-y|)}
$$ 
for all $x,y\in\bbbr^n$ and all $|\gamma|= m-1$.
\end{enumerate}
In particular, we can take $g$ such that the derivatives $D^{\gamma}g$, $|\gamma|= m-1$ are
$\lambda$-H\"older continuous simultaneously for all $\lambda\in(0,1)$.
\end{theorem}
Here $\mu(t)=O(t)$ as $t\to\infty$ means that $\mu(t)\leq Ct$ for all $t\geq t_0$.

As an application of this theorem we construct horizontal graphs in the Heisenberg group, see
Theorem~\ref{main2}. For a related construction, see also \cite{balogh}.

The paper is organized as follows. In Section~\ref{main_theorem} we prove Theorem~\ref{main}.
In Section~\ref{heisenberg} we provide a brief introduction to the Heisenberg group and then
we provide a construction of horizontal graphs based on Theorem~\ref{main}. The notation is pretty 
standard. By $C$ we will denote a general constant whose value may
change within a single string of estimates. By writing $C(n,m)$ we mean that the constant depends
on parameters $n$ and $m$ only. The symbol $C_c^m(\Omega)$ will stand for the class of compactly
supported $C^m$ functions.

\section{Proof of Theorem~\ref{main}}
\label{main_theorem}

Parts (1)-(4) of the next lemma are due to Francos \cite[Theorem~2.4]{francos}
and it is a generalization of an earlier result of Alberti \cite[Theorem~1]{alberti}.
Estimates (5) and (6) are new.

\begin{lemma}
\label{main_lemma}
Let $\Omega\subset\bbbr^n$ be open with $|\Omega|<\infty$. Let $m\geq 1$ be an integer, and let 
$f=(f_{\alpha})_{|\alpha|=m}$, be a collection of measurable functions
$f_{\alpha}:\Omega\rightarrow\bbbr$, $|\alpha|=m$. Let
$\mu:[0,\infty)\rightarrow[0,\infty)$ be a continuous function
with $\mu(0)=0$ and $\mu(t)=O(t)$ as $t\rightarrow\infty$. Let $\eps,\sigma>0$.
Then there is a function $g\in C_c^m(\Omega)$ and a compact set $K\subset\Omega$
such that
\begin{enumerate}
\item[(1)] $|\Omega\setminus K|<\eps$;
\item[(2)] $D^\alpha g(x)=f_\alpha(x)$ for all $x\in K$ and $|\alpha|=m$;
\item[(3)] 
$$
\Vert D^\alpha g\Vert_p\leq C(n,m)(\eps/|\Omega|)^{\frac{1}{p}-m}\Vert f\Vert_p
$$
for all $|\alpha|=m$ and $1\leq p\leq\infty$;
\item[(4)] $\Vert D^\gamma g\Vert_\infty<\sigma$ for all $0\leq|\gamma|<m$;
\item[(5)]
$$
|D^\gamma g(x)-D^\gamma g(y)|\leq \sigma |x-y|
$$
for all $x,y\in\bbbr^n$ and all $0\leq |\gamma|\leq m-2$;
\item[(6)]
$$
|D^\gamma g(x)-D^\gamma g(y)|\leq\frac{|x-y|}{\mu(|x-y|)}
$$ 
for all $x,y\in\bbbr^n$ and all $|\gamma|=m-1$.
\end{enumerate}
\end{lemma}
{\em Proof.}
For the proof of existence of $g\in C_c^m(\Omega)$ with properties (1)-(4), see
\cite[Theorem~2.4]{francos}. We need to prove that $g$ can be modified in such a way that (5)
and (6) are also satisfied.

Let $K'\subset\Omega$ be a compact set such that $|\Omega\setminus K'|<\eps/2$
and $f|_{K'}$ is bounded. Let $\tilde{f}=f\chi_{K'}$, where $\chi_{K'}$ is the
characteristic function of $K'$. Clearly $\Vert\tilde{f}\Vert_\infty<\infty$.
By continuity of $\mu$ we can find $\delta>0$ such that
$$
\mu(t)\leq
\frac{\eps^m}{\sqrt{n}C(n,m)|\Omega|^m\Vert\tilde{f}\Vert_\infty}
\quad
\mbox{for all $0\leq t\leq\delta$.}
$$
Here $C(n,m)$ is the constant from the inequality at (3).
In particular if $0<|x-y|\leq\delta$, then
$$
\sqrt{n}C(n,m)\eps^{-m}|\Omega|^m\Vert\tilde{f}\Vert_\infty\leq\frac{1}{\mu(|x-y|)}\, .
$$
Let
$M=\sup\{ \mu(t)/t:\, t\geq\delta\}$. $M$ is finite, because $\mu(t)=O(t)$ as $t\to\infty$.
Applying (1)-(4) to $\tilde{f}$ we can find $g\in C_c^m(\bbbr^n)$ and a compact set $K''\subset\Omega$ such that
\begin{enumerate}
\item[(1')] $|\Omega\setminus K''|<\eps/2$;
\item[(2')] $D^\alpha g(x)=f_\alpha(x)$ for all $x\in K'\cap K''$ and $|\alpha|=m$;
\item[(3')] 
$$
\Vert D^\alpha g\Vert_p \leq C(n,m) (\eps/|\Omega|)^{\frac{1}{p}-m}\Vert \tilde{f}\Vert_p
$$
for all $|\alpha|=m$ and $1\leq p\leq\infty$;
\item[(4')] 
$$
\Vert D^\gamma g\Vert_\infty < \min \left\{
\frac{\sigma}{\sqrt{n}}, \frac{1}{2M},\right\}
$$
for all $0\leq |\gamma|<m$.
\end{enumerate}

Let $K=K'\cap K''$. Then $|\Omega\setminus K|<\eps$ and it is easy to see that the function $g$
has the properties (1)-(4) from the statement of the lemma. We are left with the proof of the
properties (5) and (6). 

If $0\leq |\gamma|\leq m-2$, then (4')
yields
$$
|D^\gamma g(x)-D^\gamma g(y)|\leq
\Vert \nabla D^\gamma g\Vert_\infty |x-y|\leq \sigma|x-y|.
$$

Let now $|\gamma|=m-1$.
If $|x-y|\geq \delta$, then
$$
|D^\gamma g(x)-D^\gamma g(y)|\leq 2\Vert D^\gamma g\Vert_\infty \leq \frac{1}{M}\leq \frac{|x-y|}{\mu(|x-y|)}\, .
$$
If $0<|x-y|<\delta$, then (3') with $p=\infty$ yields
\begin{eqnarray*}
\lefteqn{|D^\gamma g(x)-D^\gamma g(y)| 
\leq
\Vert \nabla D^\gamma g\Vert_\infty |x-y|}\\
&\leq&
\sqrt{n}C(n,m)\eps^{-m}|\Omega|^m\Vert\tilde{f}\Vert_\infty |x-y|
\leq
\frac{|x-y|}{\mu(|x-y|)}\, .
\end{eqnarray*}
The proof is complete.
\hfill $\Box$

Now we can complete the proof of Theorem~\ref{main}. We follow
the argument used in \cite{francos} and \cite{moonensp},
and the only main modification is that we are using improved
estimates from Lemma~\ref{main_lemma}.

Let $U_{1}=\Omega\cap B(0,1)$. Let $V_{1}\subset\subset U_{1}$ 
be open with with
$|U_{1}\backslash V_{1}|<1/4$. Using Lemma~\ref{main_lemma}, we can
find a compact set $K_{1}\subset V_{1}$ with $|V_{1}\setminus K_{1}|<1/4$
and a function $g_{1}\in C_{c}^{m}(V_{1})$ such that
\begin{enumerate}
\item[(a)] $D^{\alpha}g_{1}(x)=f_\alpha(x)$ for all $|\alpha|=m$ and $x\in K_{1}$;
\item[(b)] 
$\left|D^{\gamma}g_{1}(x)\right|<2^{-1}\sigma\,\min{\{\rm dist}^{2}\left(x,U_{1}^{c}\right),1\}$,
for all $x\in\mathbb{R}^{n}$ and $|\gamma|<m$;
\item[(c)] 
$$
\left|D^{\gamma}g_{1}(x)-D^{\gamma}g_{1}(y)\right|\leq2^{-1}\sigma|x-y|
$$
for all $x,y\in\mathbb{R}^{n}$ and all $0\leq|\gamma|\leq m-2$;
\item[(d)]
$$
\left|D^{\gamma}g_{1}(x)-D^{\gamma}g_{1}(y)\right|\leq2^{-1}\cfrac{|x-y|}{\mu\left(|x-y|\right)}
$$
for all $x,y\in\mathbb{R}^{n}$ and all $|\gamma|=m-1$. 
\end{enumerate}

We now proceed with an inductive definition. 
Suppose that the sets $K_{1},\ldots,K_{k-1}$,
and the functions $g_{1},\ldots,g_{k-1}$ are defined, for some $k\geq 2$. Let 
$U_{k}=\Omega\cap B(0,k)\setminus(K_{1}\cup\ldots\cup K_{k-1})$.
Let $V_{k}\subset\subset U_{k}$ be open with with $|U_{k}\setminus V_{k}|<2^{-k-1}$.
Using Lemma~\ref{main_lemma}, we find a compact set $K_{k}\subset V_{k}$
with $|V_{k}\backslash K_{k}|<2^{-k-1}$and a function $g_{k}\in C_{c}^{m}(V_{k})$
such that
\begin{enumerate}
\item[(a')] $D^{\alpha}g_{k}(x)=f_\alpha(x)-\sum_{j=1}^{k-1}D^{\alpha}g_{j}(x)$, for all
$|\alpha|=m$ and $x\in K_{k}$;
\item[(b')] 
$\left|D^{\gamma}g_{k}(x)\right|<2^{-k}\sigma\,{\rm min}\left\{ {\rm dist}^{2}\left(x,U_{k}^{c}\right),1\right\}$, 
$x\in\mathbb{R}^{n}$, $|\gamma|<m$;
\item[(c')] 
$$
\left|D^{\gamma}g_{k}(x)-D^{\gamma}g_{k}(y)\right|<2^{-k}\sigma|x-y|
$$
for all $x,y\in\mathbb{R}^{n}$ and all $0\leq|\gamma|\leq m-2$;
\item[(d')]
$$
\left|D^{\gamma}g_{k}(x)-D^{\gamma}g_{k}(y)\right|\leq2^{-k}\cfrac{|x-y|}{\mu\left(|x-y|\right)}
$$
for all $x,y\in\mathbb{R}^{n}$ and all $|\gamma|=m-1$.
\end{enumerate}

We now take $g=\sum_{k=1}^{\infty}g_{k}$. We will prove that $g$ satisfies
claim of the theorem. First, to see that $g\in C^{m-1}(\mathbb{R}^{n})$, we
observe that, by (b'),
$$
\sum_{k=1}^{\infty}\left\Vert D^{\gamma}g_{k}\right\Vert _{L^{\infty}(\mathbb{R}^{n})}<\sigma
$$
for all $|\gamma|\leq m-1$, which implies $C^{m-1}$ differentiability,
and this proves (ii). Properties (iii) and (iv) now follow immediately
from (c') and (d'). Let $C=\bigcup_{k=1}^{\infty}K_{k}$. We
have $|\Omega\setminus C|=0$. 
We are left with the proof that $g$ is $m$-times differentiable at all points of $C$ and that 
$D^\alpha g=f_\alpha$, $|\alpha|=m$ on $C$.
Fix $x\in C$. Then $x\in K_{k}$ for some $k$. 
We write $g=p+q$ where $p=\sum_{j=1}^{k}g_{j}$ and
$q=\sum_{j=k+1}^{\infty}g_{j}$. Now by (a'), we have $D^{\alpha}p(x)=f_\alpha(x)$ for $|\alpha|=m$,
so we are left to show that 
$q$ is $m$-times differentiable at $x$ and $D^\alpha q(x)=0$ for $|\alpha|=m$.
Fix $|\gamma|=m-1$ and consider, for $0\neq h\in\mathbb{R}^{n}$, the difference
quotient $\left|D^{\gamma}q(x+h)-D^{\gamma}q(x)\right||h|^{-1}$.
We actually have $D^{\gamma}q(x)=0$, because $x\in K_{k}$ and 
${\rm supp}\, D^{\gamma}g_{j}\cap K_{k}=\emptyset$ for $j>k$. 
Hence
$$
\frac{|D^\gamma q(x+h)-D^\gamma q(x)|}{|h|} \leq
\frac{1}{|h|}\sum_{j=k+1}^\infty |D^\gamma g_j(x+h)|\, .
$$
If $D^\gamma g_j(x+h)\neq 0$, $j\geq k+1$, then 
$x+h\in U_j$. In this case, since also $x\in K_{k}\subset U_{j}^{c}$, we
must have ${\rm dist}(x+h,U_{j}^{c})\leq{\rm dist}(x+h,x)=|h|$. Hence
by (b') we have
$|D^\gamma g_j(x+h)|\leq 2^{-j}\sigma |h|^2$. Thus
$$
\frac{|D^\gamma q(x+h)-D^\gamma q(x)|}{|h|} \leq
\frac{1}{|h|}\sum_{j=k+1}^\infty 2^{-j}\sigma |h|^2 \leq \sigma |h|\to 0
\quad
\mbox{as $h\to 0$}
$$
which proves that the derivative $DD^\gamma q(x)=0$ equals zero
for any $|\gamma|=m-1$. This also completes the proof of (i).

In particular, if we define
$$
\mu(t)=
\begin{cases}
0 & t=0\\
|\log t|^{-1} & 0<t\leq e^{-1}\\
et & t>e^{-1},
\end{cases}
$$
then evidently $\mu$ satisfies the hypotheses above of the theorem,
and for every $\lambda\in(0,1)$ there is number $C_{\lambda}>0$ such
that
$$
\mu(t)>C_{\lambda}t^{1-\lambda},\quad t>0
$$
In that case derivatives $D^\gamma g$, $|\gamma|=m-1$ satisfy
$$
|D^\gamma g(x)-D^\gamma g(y)|\leq C_\lambda^{-1}|x-y|^\lambda
\quad
\mbox{for all $x,y\in\bbbr^n$ and $\lambda\in (0,1)$.}
$$
The proof is complete.
\hfill $\Box$

\section{The Heisenberg group}
\label{heisenberg}

In this section we will show how to use Theorem~\ref{main} to construct 
horizontal graphs in the Heisenberg group. While our construction works 
for groups $\bbbh_n$, for the sake of simplicity of notation we will restrict
to the group $\bbbh_1$; the generalization to the case of $\bbbh_n$ is
straightforward. For more information about the Heisenberg group and 
for references to results that are quoted here without proof, see
for example \cite{capogna}.

The {\em Heisenberg group}  is a Lie group
$\bbbh_1=\bbbc\times\bbbr=\bbbr^{3}$ equipped with the group law
$$
(z,t)*(z',t')=\left(z+z',t+t'+2\, {\rm Im}\,  z
  \overline{z'}\right).
$$
A basis of left invariant vector fields is given by
\begin{equation}
\label{XY}
X=\frac{\partial}{\partial x} + 2y\frac{\partial}{\partial t},\
Y=\frac{\partial}{\partial y} - 2x\frac{\partial}{\partial t},
\mbox{and}\
T=\frac{\partial}{\partial t}\, .
\end{equation}
Here and in what follows we use notation
$(z,t) = (x,y,t)$.
The Heisenberg group is equipped with the
{\em horizontal distribution} $H\bbbh_1$, which is defined at every
point $p\in\bbbh_1$ by
$$
H_p\bbbh_1={\rm span}\, \{ X(p),Y(p)\}.
$$
The distribution $H\bbbh_1$ is equipped with the left invariant metric $\boldg$
such that the vectors $X(p),Y(p)$ are
orthonormal at every point $p\in\bbbh_1$. 
An absolutely continuous curve $\gamma:[a,b]\to\bbbh_1$ is called {\em
horizontal} if $\gamma'(s)\in H_{\gamma(s)}\bbbh_1$ for almost every $s$. The
Heisenberg group $\bbbh_1$ is equipped with the {\em Carnot-Carath\'eodory
  metric} $d_{cc}$ which is defined as the infimum of the lengths of
horizontal curves connecting two given points. The lengths of curves are computed
with respect to the metric $\boldg$ on $H\bbbh_1$.
It is well known that any two points in $\bbbh_1$ can be connected by a
horizontal curve and hence $d_{cc}$ is a true metric. 
Actually, $d_{cc}$ is
topologically equivalent to the Euclidean metric. Moreover, for any compact
set $K$ there is a constant $C\geq 1$ such that
\begin{equation}
\label{SReq1}
C^{-1}|p-q|\leq d_{cc}(p,q)\leq C|p-q|^{1/2}
\end{equation}
for all $p,q\in K$. In what follows $\bbbh_1$ will be regarded as a metric space 
with metric $d_{cc}$.
The Heisenberg group is an example of a sub-Riemannian manifold \cite{gromov}.

It is often more convenient to work the {\em Kor\'anyi metric} which is
bi-Lipschitz equivalent to the Carnot-Carath\'eodory metric, but is much easier to compute.
The Kor\'anyi metric is defined by
$$
d_K(p,q)=\Vert q^{-1}*p\Vert_K,
\quad
\mbox{where}
\quad
\Vert (z,t)\Vert_K=\left(|z|^4+t^2\right)^{1/4}\, .
$$
For nonnegative functions $f$ and $g$ we write $f\approx g$ if
$C^{-1}f\leq g\leq Cf$ for some constant $C\geq 1$.
Thus bi-Lipschitz equivalence of metrics means that
$d_K\approx d_{cc}$.
A straightforward computation shows that for
$p=(z,t)=(x,y,t)$ and $q=(z',t')=(x',y',t')$ we have
\begin{equation}
\label{kor}
d_K(p,q)\approx |z-z'| +
|t-t'+2(x'y-xy')|^{1/2}.
\end{equation}

The inequality \eqref{SReq1} implies that the identity mapping
from $\bbbh_1$ to $\bbbr^{3}$ is locally Lipschitz, but its inverse is only
locally H\"older continuous with exponent $1/2$. One can prove that
the Hausdorff dimension of
any open set in $\bbbh_1$ equals $4$ and hence $\bbbh_1$ is not bi-Lipschitz homeomorphic
to $\bbbr^{3}$, not even locally.

Gromov \cite[0.5.C]{gromov} posted the following question:
{\em Given two sub-Riemannian manifolds $V$ and $W$ and $0<\alpha\leq 1$,
describe the space of
$C^\alpha$ maps $f:V\to W$.}
He also asked explicitly about the existence of H\"older continuous embeddings and 
homeomorphisms.
In particular Gromov \cite[Corollary~3.1.A]{gromov}, proved that if 
$f:\bbbr^2\to\bbbh_1$ is a $C^\alpha$-embedding,
then $\alpha\leq 2/3$ and he conjectured that $\alpha=1/2$.
This leads to a search for various surfaces in the Heisenberg group with
interesting geometric properties from the perspective of the Carnot-Carath\'eodory metric
and with suitable estimates for the H\"older continuity of a parametrization.

A problem which is related, but of independent interest, is that of finding estimates for 
the size of the characteristic set on a surface $S$ in $\bbbh_1$. We say that a point on a surface in the Heisenberg
group is {\em characteristic} if the tangent plane at this point is horizontal. 
The characteristic set $C(S)$ is the collection of all characteristic points on $S$.
In general the
Hausdorff dimension of $C(S)$ on a regular surface is small. 
Denote by $\H_E^s$ and $\dim_E$ the Hausdorff measure and the Hausdorff dimension with 
respect to the Euclidean metric.
Balogh \cite{balogh} proved that if $S$ is a $C^{2}$ surface in $\bbbh_1$, then $\dim_E(C(S))\leq 1$
and if $S$ is a $C^{1,1}$ surface, then $\dim_E C(S)<2$. On the other hand he proved that
$\bigcap_{0<\alpha<1}C^{1,\alpha}$ surfaces may satisfy $\H_E^2(C(S))>0$. 
For other related results, see \cite{franchiw}.
We should also mention the paper \cite{balogh2} that contains a construction of
horizontal fractals being graphs of $BV$ functions.

These questions motivated us in the construction of the example that we describe next (Theorem~\ref{main2}).
In what follows we will investigate surfaces in $\bbbh_1$  being graphs of continuous functions
of variables $(x,y)$. Given a function $u:\Omega\to\bbbr$, $\Omega\subset\bbbr^2$ we denote
by $\Phi(x,y)=(x,y,u(x,y))$ the canonical parametrization of the graph. We regard $\Phi$
as a mapping from $\Omega$ to $\bbbh_1$.

\begin{proposition}
\label{prop}
Suppose $\Omega\subset\bbbr^2$ is bounded. Then
$u$ is $\alpha$-H\"older continuous, $\alpha\in (0,1]$ if and only if 
$\Phi:\Omega\to\bbbh_1$ is $\alpha/2$-H\"older continuous.
\end{proposition}
{\em Proof.}
Suppose that $u$ is $\alpha$-H\"older continuous. We need to prove that
\begin{eqnarray}
\nonumber
d_K(\Phi(z),\Phi(z'))
&\approx&
|z-z'|+|u(z)-u(z')+2(x'y-xy')|^{1/2}\\
&\leq&
C|z-z'|^{\alpha/2}.
\label{a2}
\end{eqnarray}
Since $\Omega$ is bounded and $t\leq Ct^{\alpha/2}$ for $0\leq t\leq t_0$, we have
$|z-z'|\leq C|z-z'|^{\alpha/2}$ for all $z,z'\in\Omega$. Similarly boundedness of
$\Omega$ yields
\begin{eqnarray}
\nonumber
|2(x'y-xy')|^{1/2}
& \leq &
\sqrt{2}(|y|\, |x-x'|+|x|\, |y-y'|)^{1/2}
\leq
C|z-z'|^{1/2} \\
& \leq &
C'|z-z'|^{\alpha/2}\, .
\label{a3}
\end{eqnarray}
The above estimates and the $\alpha$-H\"older continuity of $u$ readily imply
\eqref{a2}.

Suppose now that $\Phi$ is $\alpha/2$-H\"older continuous i.e., \eqref{a2} is true.
The triangle inequality, \eqref{a3} and \eqref{a2} yield
\begin{eqnarray*}
|u(z)-u(z')|^{1/2}
& \leq &
|u(z)-u(z')+2(x'y-xy')|^{1/2} +
C|z-z'|^{\alpha/2} \\
& \leq &
C'|z-z'|^{\alpha/2}
\end{eqnarray*}
which in turn implies $\alpha$-H\"older continuity of $u$.
The proof is complete.
\hfill $\Box$

If $u:\Omega\to\bbbr$, $\Omega\subset\bbbr^2$ is Lipschitz continuous, then  $u$
is differentiable a.e. and hence the graph of $u$ has the tangent plane for a.e. $z\in\Omega$.
However, it cannot happen that the tangent plane to the graph is horizontal a.e.
Indeed, it is well known and easy to check that the tangent plane at
$(x,y,u(x,y))$ is horizontal if and only if
$$
\frac{\partial u}{\partial x} = 2y
\quad
\mbox{and}
\quad
\frac{\partial u}{\partial y} = -2x,
$$
but we have already checked at the beginning of this article that this system of equations
admits no Lipschitz solutions. However, we have the following result.

\begin{theorem}
\label{main2}
Let $\mu:[0,\infty)\to [0,\infty)$ be a continuous function such that
$\mu(0)=0$ and $\mu(t)=O(t)$ as $t\to\infty$. Then there is a continuous function
$u:\bbbr^2\to\bbbr$ such that
\begin{enumerate}
\item[(1)] $|u(x)-u(y)|\leq |x-y|/\mu(|x-y|)$ for all $x,y\in\bbbr^2$;
\item[(2)] $u$ is differentiable a.e.;
\item[(3)] the tangent plane to the graph of $u$ is horizontal for almost all
$(x,y)\in\bbbr^2.$
\end{enumerate}
\end{theorem}
This result is a straightforward consequence of Theorem~\ref{main} in the case of 
$m=1$ and $f_1=2y$, $f_2=-2x$. 

Note that for almost all $(x,y)\in\bbbr^n$ the
corresponding points on the surface are characteristic. 
The result is sharp -- any modulus of continuity stronger than that in (1) would
mean that the function is Lipschitz continuous and for such functions there are no
surfaces with the property (3). 
Proposition~\ref{prop} allows one to reinterpret the theorem in terms of H\"older continuous surfaces
with horizontal tangent planes. We leave details to the reader.

\end{document}